\documentclass[twoside,11pt]{article}
\usepackage{graphicx}
\usepackage{amssymb}
\usepackage{amsthm}
\newcommand{\bd}{{\rm bd}}
\newcommand{\diam}{{\rm diam}}

\newcommand{\width}{{\rm width}}
\newtheorem{thm}{Theorem}

\newtheorem{pro}{Proposition}
\font\bigbold=cmbx10 at 14 pt
\font\Bigbold=cmbx10 at 18 pt
\date{}

\title{\Bigbold 
When a spherical body of constant diameter is of constant width?}

\begin{document}

\baselineskip 17.6pt

\maketitle

\vskip-1.35cm
\centerline
{\bigbold Marek Lassak}

\pagestyle{myheadings} \markboth{Marek Lassak}{When a spherical body of constant diameter is of constant width?}

\vskip0.9cm
\noindent
{\bf Abstract.}
Let $D$ be a convex body of diameter $\delta$, where $0 < \delta < \frac{\pi}{2}$, on the $d$-dimensional sphere.
We prove that $D$ is of constant diameter $\delta$ if and only if it is of constant width $\delta$ in the following two cases.
The first case is when $D$ is smooth.
The second case is when $d=2$.



\section{Introduction}

The subject of this paper is spherical geometry (for a larger contexts see the monographs \cite{HS}, \cite{Pa} and \cite{VB}).

In the next section we recall the notion of a spherical convex body of constant width.
Shortly speaking, for a convex body $C$ on the $d$-dimensional sphere $S^d$ and any hemisphere $K$ supporting $C$ we define the {\it width of $C$ determined by $K$} as the thickness of any narrowest lune $K \cap K^*$ containing $C$.  
By a {\it body of constant width} we mean a spherical convex body whose all widths are equal.

Let $C \subset S^d$ be a convex body of diameter $\delta$.
If the spherical distance $|pq|$ of points $p, q \in C$ is $\delta$, we call $pq$ a {\it diametral chord} of $C$ and we say that $p, q$ are {\it diametrically opposed points of} $C$.
Clearly, $p, q \in \bd (C)$.
After Part 4 of \cite{LaMu-AEQ} we say that a convex body $D \subset S^d$ of diameter $\delta$ is {\it of constant diameter} $\delta$ provided for every point $p \in \bd (D)$ there exists at least one point $p' \in \bd (D)$ such that $|pp'| = \delta$ (in other words, that $pp'$ is a diametral chord of $D$).
For the known analogous notion in $E^d$ see \cite{CS}. 

Recall that in \cite {LaMu-AEQ} it is proved that a convex body on $S^d$ is of diameter $\delta \geq  \frac{\pi}{2}$ is of constant diameter if and only if it is a body of constant width $\delta$.
Moreover, there is observed that the ``if'' part holds also for $\delta < \frac{\pi}{2}$, and the problem is put {\it if every spherical body of constant diameter $\delta < \frac{\pi}{2}$ on $S^d$ is a body of constant width $\delta$?}

Our aim is to prove that every smooth spherical convex body of constant diameter $\delta$ is of constant width $\delta$, 
and that a body on the two-dimensional sphere is of constant diameter $\delta$ if and only it is of constant width $\delta$.
As a consequence of these facts, the above problem remains now open only for non-smooth bodies of constant diameter below $\frac{\pi}{2}$.

By the way, in \cite {HN1}, \cite {HN2} and \cite{L-arXiv} spherical bodies of constant width and constant diameter $\frac{\pi}{2}$ are applied for recognizing if a Wullf shape is self-dual.

\section{On spherical geometry}

By $S^d$ denote the unit sphere in the $(d+1)$-dimensional Euclidean space $E^{d+1}$, where $d\geq 2$. 
The intersection of $S^d$ with any $(k+1)$-dimensional Euclidean space, where $0 \leq k \leq d-1$, is called a {\it $k$-dimensional subsphere of $S^d$}.
For $k=1$ we call it a {\it great circle}, and for $k=0$ a {\it pair of antipodes}.
If different points $a, b \in S^d$ are not antipodes, by the {\it arc} $ab$ connecting them we mean this part of the great circle containing $a$ and $b$, which does not contain any pair of antipodes. 
By the {\it spherical distance} $|ab|$, or shortly {\it distance}, of these points we understand the length of the arc connecting them. 

By a {\it $d$-dimensional spherical ball of radius $\rho \in (0, {\frac{\pi}{2}}]$}, or shorter {\it a ball}, we mean the set of points of $S^d$ which are at the distance at most $\rho$ from a fixed point, called the {\it center} of this ball. 
For $d=2$ it is called a {\it disk}, and its boundary is called a {\it circle of radius $\delta$}.
Spherical balls of radius $\frac{\pi}{2}$ are called {\it hemispheres}.
In other words, by a {\it hemisphere} of $S^d$ we mean the common part of $S^d$ with any closed half-space of $E^{d+1}$.
We denote by $H(c)$ the hemisphere whose center is $c$.
Two hemispheres whose centers are antipodes are called {\it opposite hemispheres}.

By a {\it spherical $(d-1)$-dimensional ball of radius $\rho \in (0, {\frac{\pi}{2}}]$} we mean the set of points of a $(d-1)$-dimensional great sphere of $S^d$ at the distance at most $\rho$ from a point, called the {\it center} of this ball.  
The $(d-1)$-dimensional balls of radius $\frac{\pi}{2}$ are called {\it $(d-1)$-dimensional hemispheres}.

Let a set $C \subset S^d$ does not contain any pair of antipodes.
We say that $C$ is {\it convex} if together with every two its points, $C$ contains the whole arc connecting them. 
If the interior of a closed convex set $C$ is non-empty, we call $C$ a {\it convex body}.
Its boundary is denoted by $\bd (C)$.

If a hemisphere $H$ contains a convex body $C$ and if $p \in \bd (H) \cap C$, we say that $H$ {\it supports $C$ at $p$} or that $H$ is a {\it supporting hemisphere of $C$ at $p$}.
If exactly one hemisphere supports a convex body $C$ at its boundary point $p$, we say that $p$ is a {\it smooth point of} $\bd (C)$, and in the opposite case we say that $p$ is an {\it acute point of} $\bd (C)$ .
If every boundary point of $C \subset S^d$ is smooth, then $C$ is called {\it smooth}.
We call $C$ {\it strictly convex} if $\bd (C)$ does not contain any arc.
 
If hemispheres $G$ and $H$ of $S^d$ are different and not opposite, then $L = G \cap H$ is called {\it a lune} of $S^d$. 
This notion is considered in many books and papers.  
The parts of $\bd (G)$ and $\bd (H)$ contained in $G\cap H$ are denoted by $G/H$ and $H/G$, respectively
By the {\it thickness $\Delta (L)$ of the lune} $L = G \cap H \subset S^d$ we mean the spherical distance of the centers of the $(d-1)$-dimensional hemispheres $G/H$ and $H/G$. 

For any convex body $C \subset S^d$ and any hemisphere $K$ supporting $C$ we define the {\it width $\width_K (C)$ of $C$ determined by $K$} as the thickness of any narrowest lune $K \cap K^*$ containing $C$ (so that no lune of the form $K \cap K'$ with a smaller thickness contains $C$). 
By the {\it thickness} $\Delta (C)$ of $C$ we mean the minimum  of ${\rm width}_K (C)$ over all hemispheres $K$ supporting $C$.
Clearly, $\Delta (C)$ is nothing else but the thickness of a ``narrowest'' lune containing $C$. 
We say that $C$ is of {\it constant width} $w$ if all its widths are $w$.

The above notions are given and a few properties of lunes and convex bodies in $S^d$ are presented in \cite{L-AEQ} and \cite{LaMu-AEQ}.

\vskip0.35cm
\noindent
{\bf Lemma.}
{\it 
Let $K$ be a hemisphere of $S^d$ and let $p \in \bd (K)$.
Moreover, let $pq \subset K$ be an arc orthogonal to $\bd (K)$ with $q$ in the interior of $K$ and $|pq| < \frac{\pi}{2}$. 
Then from amongst all the lunes of the form $K \cap M$, with $q$ in the boundary of the hemisphere $M$, the lune $K \cap K_\dashv$ such that $pq$ is orthogonal to $\bd (K_\dashv)$ at $q$ has the smallest thickness.}

\vskip0.35cm
An easy proof is left to the reader.

\section{Spherical bodies of constant diameter}

The notion of a spherical body of constant diameter is recalled in the Introduction.
In this section we present a few propositions on bodies of constant diameter.

\begin{pro} \label{strictly} 
Every convex body $D \subset S^d$ of constant diameter $\delta < \frac{\pi}{2}$ is strictly convex. 
\end{pro}

\begin{proof}
Assume the opposite that $D$ is not strictly convex. 
Then $\bd (D)$ contains an arc $xz$.
Denote by $y$ its midpoint. 
Clearly, $y \in \bd (D)$.
Since $D$ is of constant diameter $\delta$, there is a point $y' \in \bd (D)$ such that $|yy'| = \delta$. 
This, $y \in xz$ and $\delta < \frac{\pi}{2}$ imply that $|y'x| > \delta$ or $|y'z| > \delta$.
Thus $\diam (D) > \delta$ in contradiction to the fact that $D$ is of constant diameter $\delta$.
Consequently, $D$ is strictly convex.
\end{proof}

\begin{pro} \label{intersect} 
Let $D \subset S^2$ be  a body of constant diameter. 
 Then every two diametral chords of $D$ intersect. 
\end{pro}

\begin{proof}
Denote the diameter of $D$ by $\delta$. 
Suppose that some diametral chords $ab$ and $cd$ of $D$ do not intersect (let for instance $a, b, d, c$ be in this order on $\bd (D)$). 
Then $abdc$ is a convex spherical non-degenerate quadrangle.
Hence $ad$ and $bc$ intersect at exactly one point.
Denote it by $x$.
Since $|ab|= \delta$ and $|cd|= \delta$, by the triangle inequality we get $|ax| + |xb| \geq \delta$ and $|cx| + |xd| \geq \delta$.
What is more, since $x \not \in ab$ and $x \not \in cd$, we get $|ax| + |xb| > \delta$ and $|cx| + |xd| > \delta$.
This leads to $|ax| + |xb| + |cx| + |xd| > 2\delta$.
So  $|ax| + |xd| > \delta$ or  $|cx| + |xb| > \delta$.
In other words, $|ad| > \delta$ or $|bc| > \delta$ in contradiction to $\diam (D) = \delta$.
Consequently, every two diametral chords of $D$ intersect.
\end{proof}

\begin{pro} \label{width}
If a hemisphere $K$ supports a convex body $D \subset S^d$ of constant diameter $\delta < \frac{\pi}{2}$ at a smooth point of its boundary, then $\width_K (D) = \delta$.
\end{pro}

\begin{proof}
Let $p$ be a smooth point of $\bd (D)$.
Since $D$ is of constant diameter, there exists a diametrally opposed point $p' \in \bd (D)$.
Hence the ball $B$ of radius $\delta$ and center $p'$ contains $D$.
Clearly, $p \in \bd (B)$.

Since $p$ is a smooth point of $\bd (D)$, we conclude that $K$ supports $B$ at $p$.

From $\diam (D) = \delta$ and $p' \in \bd(D)$ we see that the spherical ball $B'$ of radius $\delta$ centered at $p$ contains $D$ and supports it at $p'$.

Denote by $K'$ the hemisphere supporting $B'$ at $p'$. 
Clearly, we have $D \subset B \cap B' \subset K \cap K'$.

\vskip 0.75cm 

\begin{center}

\includegraphics[width=3in]{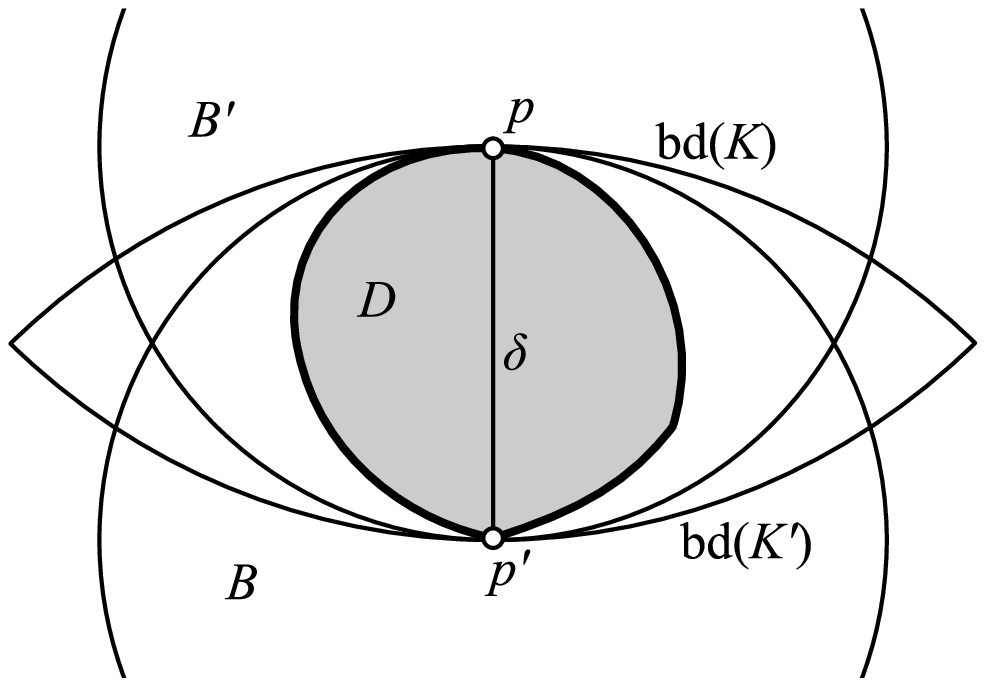} \\  

\vskip 0.25cm
{FIGURE 1. Illustration to the proof of Proposition \ref{width}} 

\end{center}
\vskip 0.1cm

Since our chord $pp'$ of length $\delta$ is orthogonal to $\bd (K)$ at $p$ and to $\bd (K')$ at $p'$, by the definition of the thickness of a lune we see that the lune $K \cap K'$ has thickness $|pp'|= \delta$. 
Consequently, by $D \subset K \cap K'$ we conclude that $\width_K (D) \leq \delta$.

By Lemma the lune $K \cap K'$ is the narrowest lune from the family ${\cal L}_1$ of lunes of the form $K \cap M$ with $p' \in \bd (M)$ containing $pp'$.  
From $pp' \subset D \subset K \cap K'$ we conclude that $K \cap K'$ is also the narrowest lune from the family ${\cal L}_2$ of all lunes the form $K \cap M$ containing $D$. 
Since every lune from ${\cal L}_2$ contains a lune from ${\cal L}_1$, every lune from ${\cal L}_2$ is of thickness at least $\delta$.
Consequently, $\width_K (D) \geq \delta$.

From the above two paragraphs we conclude that  $\width_K (D) = \delta$.  
\end{proof}

\section{Two cases in which a spherical convex body of constant diameter is of constant width}

Since the question is answered for $\delta \geq \frac{\pi}{2}$ and since, as mentioned in the Introduction, every spherical body of constant width is of constant diameter, now we concentrate on checking when a spherical body of constant diameter $\delta < \frac{\pi}{2}$ is of constant width.
The following theorem gives a partial answer.
It results immediately from Proposition \ref{width} and from the fact that every body of constant width $\delta$ is of constant diameter~$\delta$.

\begin{thm} \label{smooth}
Let $0 < \delta < \frac{\pi}{2}$. 
A spherical smooth convex body on $S^d$ is of constant diameter $\delta$ if and only if it is of constant width $\delta$. 
\end{thm}

Below is our main theorem.
Since in its proof we apply polar sets, let us recall this notion.
For a convex body $C \subset S^d$ by its {\it polar} we mean the set $C^\circ = \{r: \ C \subset H(r)\}$. 
It is easy to show that $C^\circ$ is a convex body.
Recall that $\bd (C^\circ)$ is the set of points $r$ such that $H(r)$ is a supporting hemisphere of $C$.

\begin{thm} \label{constant}
Let $0 < \delta < \frac{\pi}{2}$.
 A convex body on the two-dimensional sphere is of constant diameter $\delta$ if and only it is of constant width $\delta$.
\end{thm}

\begin{proof}
In \cite {LaMu-AEQ} there is observed that every body of constant width on $S^d$ is a body of constant diameter.

It remains to show that every body $D \subset S^2$ of constant diameter $\delta < \frac{\pi}{2}$ is of constant width $\delta$, i.e., that $\width_K (D) = \delta$ for every supporting hemisphere $K$ of $D$.

By Proposition \ref{strictly} there is exactly one point $p$ of support of $D$ by $K$.

\smallskip
When $p$ is a smooth point of $\bd (D)$, then we apply Proposition \ref{width}.

Consider the case when $p$ is an acute point of $\bd (D)$.
Denote by $H_1=H(r_1)$ the first supporting hemisphere and by $H_2 = H(r_2)$ the last supporting hemisphere of $D$ at $p$, as we go counterclockwise with the center $r$ of $H(r)$ on $\bd (D^\circ)$, see Fig 2.

For $i=1,2$ provide the arc $pp_i \subset H_i$ of length $\delta$ orthogonal to $\bd (H_i)$ at $p$.
So $p'_i \in pr_i$ for $i=1,2$.
By $B_i$ denote the ball of radius $\delta$ centered at $p_i$, where $i=1,2$.

Of course, $D \subset B_1 \cap B_2 \subset H_1 \cap H_2$.

\vskip 0.65cm 

\begin{center}

\includegraphics[width=3in]{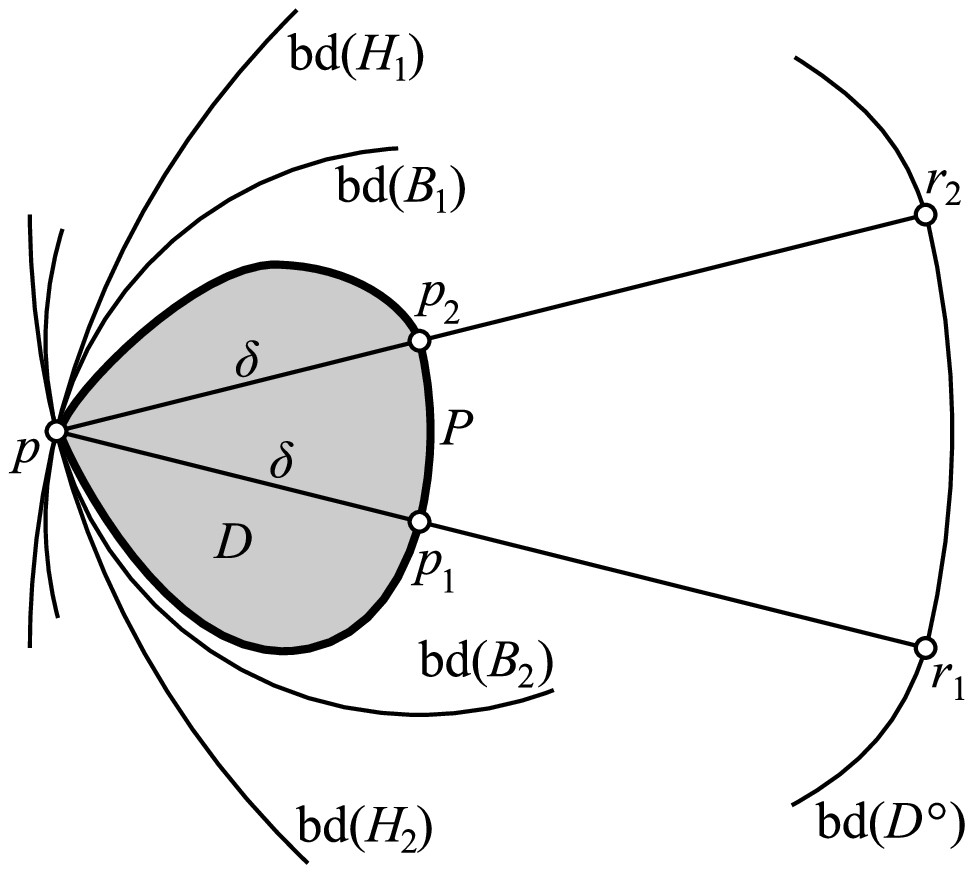} \\ 

\vskip 0.15cm
{FIGURE 2. Illustration to the proof of Theorem \ref{constant}} 

\end{center}

\vskip 0.35cm

Consider the piece $P$ of the circle with center $p$ and radius $\delta$ from $p_1$ to $p_2$.
Denote by $s_i$ the point of $\bd (D)$ in the arc $pp_i$, where $i=1,2$, and by $S$ the piece of $\bd (D)$, when going from $s_1$ counterclockwise to $s_2$.
Take any $s \in S$.
There exists a diametrally opposed point $s'$ of $D$ in $\bd (D)$. 
By Proposition \ref{intersect} the chord $ss'$ intersects the chords $pp_1$ and $pp_2$.
Hence $s'=p$.
This and $|ss'| = \delta$ imply that $s$ belongs to the circle $F$ with center $p$ and radius $\delta$. 

What is more, $s \in P$. 
The reason is that if $s \in F \setminus P$, then the disk of center $s$ and radius $\delta$ contains $D$.
Consequently, the hemisphere supporting the ball bounded by $F$ at $p$ is a supporting hemisphere of $D$ at $p$.
This contradicts the fact that  $H_1$ and $H_2$ are the first and the last supporting hemispheres of $D$ at $p$.

Consequently, for every hemisphere $H$ supporting $D$ at $p$ the chord of $D$ orthogonal to $H$ at $p$ is a diametral chord of $D$.
Thus by Lemma we get $\width_H (D) = \delta$.

In particular,  $\width_K (D) = \delta$, as required.
\end{proof}

By the proof of Theorem \ref{constant} and also Proposition \ref{width} any supporting hemisphere $H = H(r)$ of $D$ determines a unique diametral chord $pp'$ of $D$; it is orthogonal to the great circle bounding $H$. 
Moreover, the center $r$ of $H$ belongs to $\bd (D^\circ)$.
On the other hand, take any $r \in \bd (D^\circ)$.
Then $H(r)$ supports $D$ at exactly one point $p$.
Clearly, $D \subset H(r)$ . 
Consequently, for any body $D \subset S^2$ we have the one-to-one correspondence between the following objects:

- the supporting hemispheres of $D$, 

- the points of $\bd (D^\circ)$,

- the diametral chords of $D$. 
 
Every position of diametral chord is determined by a supporting hemisphere of $D$.
It is also determined by a point of $\bd (D^\circ)$; still such a point $r$ determines exactly one point $p$, and the considered diametral chord is in the arc $pp'$. 
\vskip0.15cm
A forthcoming paper is devoted to show the analogous facts as in Propositions 1--3 and Theorems 1--2 in $E^d$ (they are true also in the hyperbolic space).
The proofs of these theorems apply the parallelism, which cannot be used here for spherical bodies.

\vskip0.15cm
\noindent
Marek Lassak

\noindent
University of Science and Technology

\noindent
al. Kaliskiego 7, Bydgoszcz 85-796

\medskip
\noindent
e-mail: lassak@utp.edu.pl

\end{document}